\documentclass[12pt,leqno]{amsart}
\usepackage{amsmath,amsfonts,amssymb,amsthm,amscd,mathrsfs}
\usepackage{graphicx}
\graphicspath{ }
\usepackage{wrapfig}
\usepackage{subcaption}
\usepackage{calligra}

\numberwithin{figure}{section}

\theoremstyle{plain}
\newtheorem{theorem}{Theorem}[section]

\newtheorem{prop}[theorem]{Proposition}
\newtheorem{cor}{Corollary}[theorem]
\theoremstyle{definition}

\theoremstyle{remark}

\usepackage{mathtools}
\title[Ricci-Yamabe solitons]{Scalar curvature estimation of Generalized Ricci-Yamabe solitons}
\author[A. A. Shaikh, P. Mandal, C. K. Mondal]{Absos Ali Shaikh$^{1*}$, Prosenjit Mandal$^2$ and Chandan Kumar Mondal$^{3}$}

\address{\noindent\newline $^{1,2}$Department of Mathematics,\newline The University of Burdwan, Golapbag,\newline Burdwan-713104,\newline West Bengal, India}

\address{\noindent\newline $^3$School of Sciences,\newline Netaji Subhas Open University,\newline Durgapur Regional Center, Durgapur-713214\newline Paschim Bardhaman,\newline West Bengal, India}

\email{$^1$aask2003@yahoo.co.in, aashaikh@math.buruniv.ac.in}
\email{$^2$prosenjitmandal235@gmail.com}
\email{$^3$chan.alge@gmail.com, chandanmondal@wbnsou.ac.in}

\begin{document}
\begin{abstract}
This paper is concerned with the study of generalized gradient Ricci-Yamabe solitons. We characterize the compact generalized gradient Ricci-Yamabe soliton and find certain conditions under which the scalar curvature becomes constant. The estimation of Ricci curvature is deduced and also an isometry theorem is found in gradient Ricci-Yamabe soliton satisfying a finite weighted Dirichlet integral. Further, it is proved that a Ricci-Yamabe soliton reduces to an Einstein manifold when the potential vector field becomes concircular. Moreover, the eigenvalue and the corresponding eigenspace of the Ricci operator are also discussed in case of a Ricci-Yamabe soliton with concircular potential vector field.
\end{abstract}
\noindent\footnotetext{$^*$ Corresponding author.\\ $\mathbf{2020}$\hspace{5pt}Mathematics\; Subject\; Classification: 53C21; 53C25; 53E20.\\ 
{Key words and phrases: Riemannian manifold; Ricci-Yamabe soliton; scalar curvature; concircular vector field; potential function.} }
\maketitle
\section{Introduction and preliminaries}
In 2019, the notion of Ricci-Yamabe flow is introduced by G\"{u}ler and Cr\^{a}\c{s}m\u{a}reanu (\cite{GC2019}) by considering a scalar combination of Ricci flow and Yamabe flow originated by Hamilton (\cite{HA82, HA88}). In recent years, much effort has been devoted to the classification of self-similar
solutions of geometric flows. Ricci-Yamabe soliton is a special solution of the Ricci-Yamabe flow generated by one parameter group of diffeomorphisms and their scalings. Let $M$ be a Riemannian manifold of dimension $n(\geq 3)$ equipped with the metric $g$, then $(M,g)$ is called a Ricci-Yamabe soliton (briefly, RYS) if it preserves a smooth vector field $X$ and the Ricci curvature $Ric$ satisfying the following equation (\cite{Siddiqi20})
\begin{equation}\label{rys1}
\alpha Ric+\frac{1}{2}\mathcal{L}_Xg+(\lambda -\frac{\beta}{2} R)g=0,
\end{equation}
where $R$ and $\mathcal{L}$ indicate the scalar curvature of $M$ and Lie derivative operator respectively, and $\lambda$, $\alpha$, $\beta$ $\in \mathbb{R}$.
If $f$ is a smooth real valued function on $M$ and $X=\nabla f$, then the Ricci-Yamabe soliton is said to be a gradient Ricci-Yamabe soliton (briefly, GRYS) with potential function $f$. In this case, ($\ref{rys1}$) turns into
\begin{equation}\label{grys2}
\alpha Ric+\nabla^2f+(\lambda -\frac{\beta}{2} R)g=0,
\end{equation}
$\nabla^2f$ being the Hessian of $f$. If the potential function $f$ is constant, then the soliton is called trivial, otherwise it is non-trivial. A RYS or GRYS is said to be expanding (resp., steady or shrinking) if $\lambda>0$ (resp., $\lambda=0$ or $\lambda<0$). If $\alpha\neq 0, 1$, then a RYS or GRYS is said to be proper. An $\eta$-RYS and gradient $\eta$-RYS generalize the concept of RYS and GRYS respectively, which are defined by(\cite{Siddiqi20})
\begin{equation}\label{erys3}
\alpha Ric+\frac{1}{2}\mathcal{L}_Xg+(\lambda -\frac{\beta}{2} R)g +\mu \eta\otimes \eta=0,
\end{equation}
and
\begin{equation}\label{gerys4}
\alpha Ric+\nabla^2f+(\lambda -\frac{\beta}{2} R)g +\mu \eta\otimes \eta=0,
\end{equation}
respectively, where $\eta\otimes \eta$ denotes a (0,2)-tensor field and $\mu\in \mathbb{R}$.
In the gradient $\eta$-RYS equation ($\ref{gerys4}$), if we replace $\eta\otimes \eta$ by $df\otimes df$, where $df$ is the dual $1$-form of the vector filed $\nabla f$, then ($\ref{gerys4}$) takes the following form 
\begin{equation}\label{gerys5}
\alpha Ric+\nabla^2f+(\lambda -\frac{\beta}{2} R)g +\mu df\otimes df=0.
\end{equation}
 If $(M,g)$ obeys the equation (\ref{gerys5}), then it is called generalized gradient Ricci-Yamabe solition(\cite{Dey21}).

In local coordinate system (\ref{gerys5}) can be written as
\begin{equation}\label{gerys6}
\alpha R_{ij}+\nabla_i \nabla_jf+(\lambda -\frac{\beta}{2} R)g_{ij} +\mu \nabla_i f \nabla_jf=0.
\end{equation}
We note that, for $\mu=0$, a gradient $\eta$-RYS (resp., $\eta$-RYS) is simply a GRYS (resp., RYS).

 The concept of Ricci-Yamabe soliton generalizes a huge class of soliton like equations such as: (i) for $\alpha=1$, $\beta=0$, it is a Ricci soliton (\cite{HA82,HA88, CA20, SD22}), (ii) for $\alpha=0$, $\beta=2$, it represents Yamabe soliton (\cite{HA82,HA88}), (iii) for $\alpha=1$, $\beta=2\rho$, it is $\rho$-Einstein soliton (\cite{Catino16,Catino15}), etc. Recently, the study of various types of solitons and their generalizations have been the focus of
 attention of many mathematicians. In (\cite{chen2015}), Chen classified Ricci
solitons with concircular potential vector field and deduced that the only concircularly flat manifold admitting a Ricci soliton with concircular potential field is the Euclidean space. Global problems about various types of vector fields can be found in Mike\v s et al. (\cite{metal,mhg}) and also references therein. In (\cite{Catino15}), Catino et al. established a rigidity result and also proved some fundamental identities for gradient shrinking  $\rho$-Einstein solitons. Later, Shaikh et al. (\cite{AAP2021, SMM20,SMM2022}), studied $\rho$-Einstein solitons and proved some interesting results.

 Hence inspiring by the above studies, we obtain the following results:
\begin{prop}\label{thm6}
 Let $(M,g)$ be an $n$-dimensional compact generalized gradient Ricci-Yamabe soliton. If the scalar curvature $R$ agrees with the potential function $f$ at its points of extrema and $n\beta>2\alpha$, then the following holds:\\
  (i) the scalar curvature $R$ is constant, and\\
  (ii) the soliton is trivial. 
  \end{prop}
  In view of Proposition \ref{thm6}, we have the following corollary:
  \begin{cor}\label{co3}
  Let $(M,g)$ be an $n$-dimensional compact generalized gradient Ricci-Yamabe soliton. If the scalar curvature $R$ agrees with the potential function $f$ at its points of extrema and $n\beta>2\alpha$, then the scalar curvature is negative (resp., vanishes, positive) according as the soliton is shrinking (resp., steady, expanding).
  \end{cor}
  We generalize the work of Catino et al. (\cite{Catino15}, Lemma 2.2) and prove the following
  fundamental identities:
    \begin{theorem}\label{thm1}
   In a generalized gradient Ricci-Yamabe soliton $(M,g)$ of dimension $n(\geq 3)$, the following relations holds:
   \begin{equation}\label{gerys8}
   \alpha R+\Delta f+(\lambda -\frac{\beta}{2} R)n+\mu |\nabla f|^2=0,
   \end{equation}
   \begin{equation}\label{gerys9}
   \{\alpha-\beta(n-1)\}\nabla R+2\mu\{\alpha R+(\lambda-\frac{\beta}{2}R)(n-1)\}\nabla f =2(\mu \alpha+1)Ric(\nabla f,.),
   \end{equation}
   \begin{eqnarray}
  \nonumber \{\alpha-\beta(n-1)\}\Delta R+\{2\mu \alpha -2\mu \beta (n-1)-1\}g(\nabla R,\nabla f)
   = 2\mu \{\alpha R+(n-1)(\lambda-\frac{\beta}{2}R)\}\\
   \times\{\alpha R+n(\lambda-\frac{\beta}{2}R)\}
   -2(\mu \alpha+1) \{\alpha |Ric|^2+R(\lambda-\frac{\beta}{2}R)\}.
   \end{eqnarray}
    \end{theorem}
     For Ricci-Yamabe soliton, i.e., for $\mu=0$, we deduce the following:
     \begin{cor}\label{thm7}
     If $(M,g)$ is a gradient Ricci-Yamabe soliton of dimension $n$, then the following relations hold:
     \begin{equation}\label{grys11}
     \alpha R+\Delta f+(\lambda -\frac{\beta}{2} R)n=0,
     \end{equation}
     \begin{equation}\label{grys12}
     \{\alpha-\beta(n-1)\}\nabla R =2 Ric( \nabla f,.),
     \end{equation}
     \begin{equation}\label{grys13}
    \{\alpha-\beta(n-1)\}\Delta R
     =g(\nabla R,\nabla f)-2 \{\alpha |Ric|^2+R(\lambda-\frac{\beta}{2}R)\}.
     \end{equation}
      \end{cor}
      We note that for $\alpha=1$, $\beta=2\rho$, $(M,g)$ is $\rho$-Einstein soliton, thus the above result generalizes the  Lemma 2.2 of Catino et al. (\cite{Catino15}).
 \begin{theorem}\label{thm5}
  If $(M,g)$ is an $n$-dimensional compact gradient steady Ricci-Yamabe soliton, then the scalar curvature $R$  satisfies the following integral inequality
  \begin{equation}
  k \int_{M} R^2 \geq \int_{M}Ric(\nabla f, \nabla f),
  \end{equation}
  for some scalar $k=\frac{n-1}{n} (\frac{\beta n}{2}-\alpha)^2.$ 
   \end{theorem}
\begin{theorem}\label{thm4}
 Let $(M,g)$ be an $n$-dimensional complete connected non-trivial gradient Ricci-Yamabe soliton satisfying the following Dirichlet integral condition
 \begin{equation}\label{Dirichlet condition}
 \int_{M\backslash B(q,r)}|\nabla f|^2d(x,q)^{-2}<\infty,
 \end{equation}
 $B(q,r)$ being the open ball with center at $q\in M$ and radius $r>0$ and $d(x,q)$ is the distance function from some fixed point $q\in M$. If any one of the following
 \begin{enumerate}
 \item[i)] $Ric\leq 0$ and either $\beta>\alpha$ or $\beta<\frac{\alpha}{(n-1)}$,
 \item[ii)] $Ric\geq 0$ and  $\frac{\alpha}{(n-1)}<\beta<\alpha$,
 \end{enumerate}
 holds, then $M$ is isometric to $N\times \mathbb{R}$, where $N$ is totally geodesic submanifold of $M$.
  \end{theorem}
  The Theorem \ref{thm4}, implies the following corollary:
  \begin{cor}\label{co2}
   Let $(M,g)$ be an $n$-dimensional complete non-trivial gradient Ricci-Yamabe soliton with $\alpha \neq 0$ and satisfying $(\ref{Dirichlet condition})$. If any one of the following
   \begin{enumerate}
   \item[i)] $Ric\leq 0$ and either $\beta>\alpha$ or $\beta<\frac{\alpha}{(n-1)}$,
   \item[ii)] $Ric\geq 0$ and  $\frac{\alpha}{(n-1)}<\beta<\alpha$,
   \end{enumerate}
   holds, then the soliton is steady Ricci flat.
  \end{cor}
  We use the technique of (\cite{chen2015, DC2015}) and prove the following results for Ricci-Yamabe solitons:
  \begin{theorem}\label{thm3}
   Let $(M,g)$ be an $n$-dimensional Ricci-Yamabe soliton possessing concircular potential vector field $X$ with $\alpha \neq 0$. Then $M$ is Einstein and the scalar curvature is given by $\frac{2n(\lambda + \phi)}{\beta-2\alpha}$.
    \end{theorem}
    By virtue of Theorem \ref{thm3}, the following corollary can be stated:
    \begin{cor}\label{co1}
    Let $(M,g)$ be an $n$-dimensional Ricci-Yamabe soliton possessing concircular potential vector field $X$ with $\alpha \neq 0$. Then the Ricci-Yamabe soliton is expanding, steady or shrinking according as $\phi< \frac{\beta-2\alpha}{2n}R$, $\phi=\frac{\beta-2\alpha}{2n}R$ or $\phi>\frac{\beta-2\alpha}{2n}R$.
    \end{cor}
    \begin{theorem}\label{thm2}
     Let $(M,g)$ be an $n$-dimensional Ricci-Yamabe soliton possessing concircular potential vector field $X$ with $\alpha \neq 0$. Then the space $\chi (M)$ of all smooth vector fields is the eigenspace of the Ricci operator $Q$ corresponding to the eigenvalue $\frac{\beta R-2\phi-2\lambda}{2\alpha}$.
      \end{theorem}
\section{Proof of the results}
\begin{proof}[\textbf{Proof of Proposition \ref{thm6}}]
Tracing the equation (\ref{gerys5}), we have
\begin{equation}\label{gerys10}
\Delta f+\mu |\nabla f|^2=- \alpha R-(\lambda -\frac{\beta}{2} R)n.
 \end{equation}
 At the point of maximum value of $f$, $\nabla f=0$ and $\Delta f \leq 0 $.
 Thus (\ref{gerys10}) implies that
 $$- \alpha R_{max}-(\lambda -\frac{\beta}{2} R_{max})n\leq 0,$$
 where $R_{max}$ represents the maximum value of $R $. Therefore, a simple calculation yields $$R_{max}\leq \frac{2n \lambda }{n\beta-2\alpha}.$$
 Similarly, at the point of minimum value of $f$, we have
 $$R_{min}\geq \frac{2n \lambda }{n\beta-2\alpha}.$$
 Hence, $$f_{min}=R_{min}=R_{max}=f_{max}= \frac{2n \lambda }{n\beta-2\alpha}.$$
 Therefore, $R=f=\frac{2n \lambda }{n\beta-2\alpha}$, a constant. This completes the proof of the results.
\end{proof}
\begin{proof}[\textbf{Proof of Theorem \ref{thm1}}]
Taking the trace of (\ref{gerys5}), we get the identity (\ref{gerys8}) as follows
\begin{equation}\label{eq2}
\alpha R+\Delta f+(\lambda -\frac{\beta}{2} R)n+\mu |\nabla f|^2=0.
\end{equation}
From the commutative equation for covariant derivative, we have
\begin{equation}\label{com7}
R_{ij}\nabla_jf=\Delta\nabla_if -\nabla_i\Delta f.
\end{equation}
By the contracted second Bianchi identity and equation (\ref{gerys6}), we obtain
\begin{eqnarray*}
\Delta\nabla_if=\nabla_j\nabla_j\nabla_i f&=&\nabla_j\{-(\lambda -\frac{\beta}{2} R)g_{ij}-\alpha R_{ij} -\mu \nabla_i f \nabla_jf\}\\
&=& \frac{\beta}{2}\nabla_i R-\frac{\alpha}{2} \nabla_i R-\mu (\nabla_i\nabla_j f \nabla_jf+\nabla_if \nabla_j\nabla_jf,
\end{eqnarray*}
and
\begin{eqnarray*}
\nabla_i\Delta f&=&\nabla_i\{-\alpha R-(\lambda -\frac{\beta}{2} R)n-\mu |\nabla f|^2\}\\
&=&-\alpha \nabla_iR+\frac{\beta}{2}n \nabla_i R-2\mu \nabla_jf \nabla_i \nabla_jf.
\end{eqnarray*}
$$$$
Therefore, (\ref{com7}) yields
\begin{equation}
\frac{\beta}{2}(1-n)\nabla_i R+\frac{\alpha}{2} \nabla_i R+\mu (\nabla_jf\nabla_i\nabla_j f-\nabla_if \Delta f)=R_{ij}\nabla_jf.
\end{equation}
Thus using (\ref{gerys6}) and (\ref{gerys8}), the above equation can be written as follows
\begin{equation}
\{\alpha-\beta(n-1)\}\nabla_i R+2\mu\{\alpha R+(\lambda-\frac{\beta}{2}R)(n-1)\}\nabla_if =2(\mu \alpha+1)R_{ij}\nabla_jf,
\end{equation}
which follows the identity (\ref{gerys9}). Now, taking covariant derivative $\nabla_l$ of the above equation, we get
\begin{eqnarray*}
&&\{\alpha-\beta(n-1)\}\nabla_l\nabla_i R+2\mu\{\alpha  (\nabla_l R\nabla_i f+R\nabla_l\nabla_i f)+(\lambda-\frac{\beta}{2}R)(n-1)\}\nabla_l\nabla_if\\
&&-\frac{\beta}{2}(n-1)\nabla_l R\nabla_i f
= 2(\mu \alpha+1)(\nabla_l R_{ij}\nabla_jf +R_{ij} \nabla_l\nabla_j f).
\end{eqnarray*}
Taking trace in $l$ and $i$, and using the contracted second Bianchi identity, we obtain 
\begin{eqnarray*}
&&\{\alpha-\beta(n-1)\}\Delta R+2\mu\{\alpha  g(\nabla R,\nabla f+\alpha R\Delta f)+(\lambda-\frac{\beta}{2}R)(n-1) \Delta f\\&&-\frac{\beta}{2}(n-1)g(\nabla R,\nabla f)\}
= 2(\mu \alpha+1)(\frac{1}{2}g(\nabla R,\nabla f)+R_{ij} \nabla_i\nabla_jf).
\end{eqnarray*}
Simplifying the above equation, we get
\begin{eqnarray*}
&&\{\alpha-\beta(n-1)\}\Delta R+\{2\mu \alpha -\mu \beta (n-1)-(\mu \alpha+1) \}g(\nabla R,\nabla f)\\
&&+2\mu \{\alpha R+(\lambda-\frac{\beta}{2}R)(n-1)\} \Delta f
= 2(\mu \alpha+1) R_{ij} \nabla_i\nabla_jf.
\end{eqnarray*}
Again using (\ref{gerys6}) and (\ref{gerys8}), the above equation takes the following form
\begin{eqnarray*}
&&\{\alpha-\beta(n-1)\}\Delta R+\{2\mu \alpha -\mu \beta (n-1)-(\mu \alpha+1) \}g(\nabla R,\nabla f)\\
&=& 2(\mu \alpha+1) \{-\alpha |Ric|^2-R(\lambda-\frac{\beta}{2}R)\}
+2\mu \{\alpha R+(n-1)(\lambda-\frac{\beta}{2}R)\}\{\alpha R+n(\lambda-\frac{\beta}{2}R)\}\\&&-2\mu(\mu \alpha+1) R_{ij} \nabla_i f \nabla_jf+2\mu^2 \{\alpha R+(n-1)(\lambda-\frac{\beta}{2}R)\}\nabla_i f \nabla_i f.
\end{eqnarray*}
Now, using (\ref{gerys9}), we obtain
\begin{eqnarray*}
\{\alpha-\beta(n-1)\}\Delta R+\{2\mu \alpha -2\mu \beta (n-1)-1\}g(\nabla R,\nabla f)\\
= 2\mu \{\alpha R+(n-1)(\lambda-\frac{\beta}{2}R)\}\{\alpha R+n(\lambda-\frac{\beta}{2}R)\}-2(\mu \alpha+1) \{\alpha |Ric|^2+R(\lambda-\frac{\beta}{2}R)\}.
\end{eqnarray*}
This completes the proof.
\end{proof}
\begin{proof}[\textbf{Proof of Theorem \ref{thm5}}]
The Bochner formula together with the equation (\ref{grys11}) and $\lambda=0$, yields
$$\frac{1}{2}\Delta|\nabla f|^2=|\nabla^2 f|^2+Ric(\nabla f, \nabla f)+(\frac{\beta n}{2}-\alpha)\langle \nabla R,\nabla f\rangle.$$
Now using the fact, $|\nabla^2 f|^2\geq \frac{1}{n}(\Delta f)^2$, we obtain
\begin{eqnarray*}
 \frac{1}{2}\Delta|\nabla f|^2&\geq
 & \frac{1}{n}(\frac{\beta n}{2}-\alpha)^2 R^2+Ric(\nabla f, \nabla f)+(\frac{\beta n}{2}-\alpha)\langle \nabla R,\nabla f\rangle.
 \end{eqnarray*}
 Using divergence theorem, we get
 \begin{eqnarray*}
  0&\geq
  &  \int_{M} \frac{1}{n}(\frac{\beta n}{2}-\alpha)^2 R^2+\int_{M} Ric(\nabla f, \nabla f)-(\frac{\beta n}{2}-\alpha) \int_{M} R \Delta f\\
  &=& \int_{M}(\frac{1}{n}-1)(\frac{\beta n}{2}-\alpha)^2 R^2+\int_{M} Ric(\nabla f, \nabla f).\\
  \end{eqnarray*}
  Therefore,$$\frac{n-1}{n} (\frac{\beta n}{2}-\alpha)^2 \int_{M} R^2\geq\int_{M} Ric(\nabla f, \nabla f).$$
  This follows the result.
\end{proof}
\begin{proof}[\textbf{Proof of Theorem \ref{thm4}}]

From Bochner formula and the equation (\ref{grys11}), we have
$$\frac{1}{2}\Delta|\nabla f|^2=|\nabla^2 f|^2+Ric(\nabla f, \nabla f)+\frac{\beta n-2\alpha}{2}\langle \nabla R,\nabla f\rangle.$$
Using (\ref{grys12}), the above equation yields
\begin{eqnarray}\label{1}
 \nonumber\frac{1}{2}\Delta|\nabla f|^2 &=&|\nabla^2 f|^2+ Ric(\nabla f, \nabla f)+\left\{\frac{\beta n-2\alpha}{\alpha-\beta(n-1)}\right\}Ric(\nabla f, \nabla f)\\
&=&|\nabla^2 f|^2+\left\{\frac{\beta -\alpha}{\alpha-\beta(n-1)}\right\}Ric(\nabla f, \nabla f).
\end{eqnarray}
As in \cite{Cheeger}, let $\zeta_r\in C_0^\infty(B(q,2r))$ for $r>0$, such that
\begin{eqnarray}\label{cuttoff}\left\{ \begin{array}{lllll}
0\leq\zeta_r\leq1 & {\rm in}\,\,B(q,2r)\\
\,\,\,\,\,\zeta_r=1 & {\rm in}\,\, B(q,r)\\
\,\,\,\,\,|\nabla\zeta_r|^2\leq\frac{C}{r^2} & {\rm in}\,\,B(q,2r)\\
\,\,\,\,\,\Delta\zeta_r\leq\frac{C}{r^2} & {\rm in}\,\,B(q,2r),\\
\end{array}\right.
\end{eqnarray}
where $C>0$ is a constant.
Integration by parts and our assumption together imply that
$$\int_{B(p,2r)}\zeta^2_r\Delta|\nabla f|^2=\int_{B(p,2r)}|\nabla f|^2\Delta\zeta^2_r\leq\int_{{B(p,2r)}\backslash B(q,r)}\frac{C}{r^2}|\nabla f|^2\rightarrow 0,$$
as $r\rightarrow\infty$. 
Hence, (\ref{1}) yields
\begin{equation}\label{hessian}
\int_M|\nabla^2 f|^2+\left\{\frac{\beta -\alpha}{\alpha-\beta(n-1)}\right\}\int_M Ric(\nabla f,\nabla f)=0.
\end{equation}
By the given conditions 
$$\left\{\frac{\beta -\alpha}{\alpha-\beta(n-1)}\right\} Ric(\nabla f,\nabla f)\geq 0.$$
Thus, from (\ref{hessian}) it follows that
$\nabla^2 f$ vanishes. Hence, (see, \cite{Sakai96}, Lemma 2.3) we conclude the following:\\
(i) the smooth function $f$ is affine and\\
(ii) $\nabla f$ is a Killing vector field with $|\nabla f|$= constant.  
Since $M$ admits a non constant affine function $f$, it is isometric to $N\times \mathbb{R}$ (see, \cite{Innami82}), where $N$ is totally geodesic submanifold of $M$.
\end{proof}
\begin{proof}[\textbf{Proof of Corollary \ref{co2}}]
By the given conditions from equation (\ref{hessian}), we obtain $Ric(\nabla f,\nabla f)=0$ and $\nabla^2 f=0$. Using these in the gradient Ricci-Yamabe soliton equation (\ref{grys2}), we get $\lambda=\frac{\beta}{2}R$, which follows that $Ric \equiv 0$ and $\lambda=0$.
\end{proof}
\begin{proof}[\textbf{Proof of Theorem \ref{thm3}}]
 Since $X$ is concircular, there exists a smooth function $\phi$ on $M$ such that $\nabla_Y X=\phi Y$, for all $Y \in \chi(M)$. Then we have $$\frac{1}{2}(\mathcal{L}_Xg) (Y,Z)=\phi g( Y,Z).$$
 Thus from equation $(\ref{rys1})$, we obtain
\begin{equation}\label{concircular2}
\alpha Ric(Y,Z)+\phi g( Y,Z)+(\lambda -\frac{\beta}{2} R)g(Y,Z)=0,
\end{equation}
which gives
\begin{equation}\label{cone1}
 Ric(Y,Z)=\frac{\beta R-2\phi-2\lambda}{2\alpha} g(Y,Z).
\end{equation}
Since $n\geq 3$, using the contracted second Bianchi identity, we conclude that $M$ is Einstein. Thus $\frac{\beta R-2\phi-2\lambda}{2\alpha}$ is constant. Hence from (\ref{cone1}), we obtain that the constant scalar curvature $R$ is given by $\frac{2n(\lambda + \phi)}{\beta-2\alpha}$.
\end{proof}
\begin{proof}[\textbf{Proof of Theorem \ref{thm2}}]
 Since $X$ is concircular, there exists a smooth function $\phi$ on $M$ such that $\nabla_Y X=\phi Y$, for all $Y \in \chi(M)$. Then we have $$\frac{1}{2}(\mathcal{L}_Xg) (Y,Z)=g(\phi Y,Z).$$
 Thus from equation $(\ref{rys1})$, we obtain
\begin{equation*}
\alpha Ric(Y,Z)+g(\phi Y,Z)+(\lambda -\frac{\beta}{2} R)g(Y,Z)=0.
\end{equation*}
Then using the definition of symmetric (1,1)-Ricci operator $Q$, i.e., $Ric(Y,Z)=g((Q)Y,Z)$, we get
\begin{equation*}
\alpha g((Q)Y,Z)+g(\phi Y,Z)+(\lambda -\frac{\beta}{2} R)g(Y,Z)=0,
\end{equation*}
which yields
\begin{equation}\label{ro1}
\alpha Q(Y)+\phi Y+(\lambda -\frac{\beta}{2} R)Y=0.
\end{equation}
Now, from (\ref{ro1}), we obtain
\begin{equation}\label{ro2}
 Q(Y)=\frac{\beta R-2\phi-2\lambda}{2\alpha} Y.
\end{equation}
As $Y$ is arbitrary in $\chi(M)$, from Theorem \ref{thm2} the result follows.
\end{proof}
\section*{Acknowledgment}
 The second author gratefully acknowledges to the CSIR(File No.:09/025(0282)/2019-EMR-I), Govt. of India for the award of Junior Research Fellow.
\section{Data Availability Statement}
Our manuscript has no associated data.
\section{Conflict of Interest}
The authors declare no conflict of interest.

\end{document}